\documentclass[conference] {IEEEtran}
\usepackage{graphicx} 
\usepackage{amsmath}
\usepackage{amsfonts}
\usepackage{amsmath,latexsym,amsfonts,enumerate}
\usepackage{amssymb,amsthm}
\usepackage{mathtools}
\usepackage{csquotes}
\usepackage{tikz-cd}
\usepackage{tikz}
\usepackage{comment}
\usepackage[T1]{fontenc}
\usepackage{cleveref}
\usepackage{lmodern}
\usepackage[total={6in,10in},left=1.5in,top=0.5in,includehead,includefoot]{geometry}
\usepackage{microtype}
\usepackage[nodisplayskipstretch]{setspace}
\usepackage{amsmath}
\usepackage{amsfonts}
\usepackage{mathtools}
\usepackage{tensor}
\everymath{\displaystyle}

\def\R{\mathbb{R}}

\DeclareMathOperator{\dispr}{{\mathbf{d}_{PR}}}
\DeclareMathOperator{\disper}{{\mathbf{d}}}
\DeclareMathOperator{\lp}{\langle}

\DeclareMathOperator{\rp}{\rangle}

\newtheorem{theor}{Theorem}

\newtheorem{cor}[theor]{Corollary}
\newtheorem{defin}[theor]{Definition}
\newtheorem{rem}[theor]{Remark}

\newtheorem*{ques*}{Question}
\newtheorem*{theor*}{Theorem}
\newtheorem*{lem*}{Lemma}
\newtheorem*{prop*}{Proposition}
\newtheorem*{cor*}{Corollary}
\newtheorem*{defin*}{Definition}
\newtheorem*{rem*}{Remark}
\newtheorem*{rems*}{Remarks}
\newtheorem*{ex*}{Example}
\newtheorem*{notat*}{Notations}
\newtheorem*{claim*}{Claim}

\title{Relationships between the Phase Retrieval Problem and Permutation Invariant Embeddings}
\author{}
\date{March 31, 2023}

\newcommand{\ignore}[1]{}

\newcommand{\ip}[2]{{\lp#1,#2\rp}}
\newcommand{\norm}[1]{{\|#1\|}}

\begin{document}

\author{
\IEEEauthorblockN{Radu Balan and Efstratios Tsoukanis}
\IEEEauthorblockA{ \\ Department of Mathematics, 
University of Maryland, 
College Park, MD 20742\\
Emails: rvbalan@umd.edu , etsoukan@umd.edu
}
}

\maketitle

\begin{abstract}
This paper discusses the connection between the phase retrieval problem and permutation invariant embeddings. We show that the real phase retrieval problem for $\R^d/O(1)$ is equivalent to Euclidean embeddings of the quotient space $\R^{2\times d}/S_2$ performed by the sorting encoder introduced in an earlier work. In addition, this relationship provides us with inversion algorithms of the orbits induced by the group of permutation matrices.  
\end{abstract}

\section{Introduction}

The phase retrieval problem has a long and illustrious history involving several Nobel prizes along the way. 
The issue of reconstruction from magnitude of frame coefficients 
is related to a significant number of problems that appear in
separate areas of science and engineering. Here is an incomplete list of some of these applications and reference papers:  crystallography \cite{KaKa,GS,fienup}; 
ptychography \cite{ZhengBook,saab20a}; 
source separation and inverse problems
\cite{OLKP79,CISS10}; optical data processing \cite{winkler09}; mutually
unbiased bases \cite{DGS77,SPIE2007a}, quantum state tomography \cite{HMW11,KW15}; 
low-rank matrix completion problem \cite{CSV12,CESV12}; tensor algebra and systems of multivariate 
polynomial equations \cite{LV12,CL12,CEHV13}; signal generating models \cite{shamma2004,Stoica07},  bandlimited
functions \cite{thakur10,CCSW16}, radar ambiguity problem \cite{Jaming,Jaming2010}, learning and scattering networks \cite{M12,BM14,Romberg16a}.

In \cite{balan06}, this problem was shown to be a special form of the following setup.  Let $H$ denote a real or complex vector space and let $A=\{a_i\}_{i\in I}$ be a frame for $H$. The phase retrieval problem asks whether the map $H\ni x\mapsto \alpha_{A}(x)=\{|\ip{x}{a_i}|\}_{i\in I}\in l^2(I)$ determines $x$ uniquely up to a unimodular scalar.

In this paper we focus on the finite dimensional real case of this problem (see also \cite{bcmn}), namely when $H=\R^{d}$.
In this case, a frame $\mathcal{A}=\{a_1,\ldots,a_D\}\subset \R^d$ is simply a spanning set. The group $O(1)=\{-1,+1\}$ acts on $H$ by scalar multiplication. Let $\hat{H}=H/O(1)$ denote the quotient space induced by this action, where the equivalence classes (orbits) are
\[ [x]=\{x,-x\}  ~~,~~ {\text{for}} ~x\neq 0~~,~~[x]=\{0\}~~,~~\text{for}~x=0.\]
The analysis operator for this frame is 
\begin{equation}
    \label{eq1}
T_A:H\rightarrow\R^D ~~,~~ T_A(x) = (\ip{x}{a_k})_{k=1}^D.
\end{equation}
The relevant nonlinear map $\alpha_A$ is given by taking the absolute value of entries of $T_A$:
\begin{equation}
    \label{eq2}
\alpha_A: H\rightarrow \R^{D} ~~,~~
\alpha_A(x)=(|\ip{x}{a_k}|)_{ k=1}^{ D}.
\end{equation} 

Notice $\alpha_A$ produces a well-defined map on $\hat{H}$, which, with a slight abuse, but for simplicity of notation, will be denoted also by $\alpha_A$. Thus $\alpha_A([x])=\alpha_A(x)$.

Another customary notation that is often employed: a frame is given either as an indexed set of vectors, $\mathcal{A}=\{a_1,\ldots,a_D\}$, or through the columns of a $d\times D$ matrix $A$. 
The matrix notation is not canonical, but this is not an issue here. We always identify $H=\R^d$ with its columns vector representation in its canonical basis.
\begin{defin}
We say  that (the columns of a matrix) $A \in \R^{d \times D}$ form/is a {\em phase retrievable frame}, if  $\alpha_A : \widehat{\R}^d \to \R^D$, $\alpha_A(x)= (|\ip{x}{a_k}|)_{k=1}^D$ is an injective map (on the quotient space).
\end{defin}

In a different line of works \cite{Cahill19,balan2022permutation,cahill22,dym22}
 it was recognized that the phase retrieval problem is a special case of Euclidean representations of metric spaces of orbits defined by certain unitary group actions on Hilbert spaces. Specifically, the setup is as follows. Let $V$ denote a Hilbert space, and let $G$ be a group acting unitarily on $V$. Let $\hat{V}=V/G$ denote the metric space of orbits, where the quotient space is induced by the equivalence relation $x,y\in V$, $x\sim y$ iff $y=g.x$, for some $g\in G$. Here $g.x$ represents the action of the group element $g\in G$ on vector $x$. For the purposes of this paper we specialize to the finite dimensional real case, $V=\R^{n\times d}$ and $G=S_n$, is the group of $n\times n$ permutation matrices acting on $V$ by {\em left multiplication}.
 Other cases are discussed in aforementioned papers. In particular, in \cite{balan2022permutation} the authors have shown a deep connection to graph deep learning problems. In \cite{cahill22}, the authors linked this framework to certain graph matching problems and more.
The bi-Lipschitz Euclidean embedding problem for the finite dimensional case is as follows. Given $\hat{V}=V/G$, construct a map $\beta:V\rightarrow\R^m$ so that, (i) $\beta(g.x)=\beta(x)$ for all $g\in G$, $x\in V$, and (ii) for some $0<A\leq B<\infty$, and for all $x,y\in V$,
\begin{equation} 
\label{eq3}
A\, \disper([x],[y]) \leq \norm{\beta(x)-\beta(y) } \leq B\,  \disper([x],[y]) \end{equation}
where $\disper([x],[y])=\inf_{g\in G}\norm{x-g.y}_V$ is the {\em natural metric} on the quotient space $\hat{V}$.

 In \cite{balan2022permutation} the following embedding was introduced. Let $A\in\R^{d\times D}$ 
 be a fixed matrix (termed as {\em key}) whose columns are denoted by $a_1,\ldots,a_D$. The induced encoder $\beta_A:V\rightarrow \R^{n\times D}$ is defined by
 \begin{equation} \label{eq4}
\beta_A(X) = \downarrow(XA) = \left[ 
\begin{array}{ccc}
     \Pi_1 X a_1 &  \cdots & \Pi_D X a_D  
\end{array}
\right]
 \end{equation}
where $\Pi_k\in S_n$ is the permutation matrix that sorts in decreasing order the vector $Xa_k$.
It was shown in \cite{balan2022permutation} that, for $D$ large enough, $\beta_A$ provides a bi-Lipschitz Euclidean embedding of $\hat{V}$. This motivates the following definition.

\begin{defin}
We say that $A\in\R^{d\times D}$ is a {\em universal key} for $\R^{n\times d}$ if $\beta_A:\widehat{\R^{n\times d}}\rightarrow\R^{n\times D}$, $\beta_A(X)\!=\downarrow\!(XA)
$ is an injective map (on the quotient space).
\end{defin}

The purpose of this paper is to show the equivalence between the real phase retrieval problem, specifically the embedding $\alpha_A$, and the
permutation invariant embedding $\beta_A$ defined above, in the special case $n=2$. 

\section{Main Results}
Recall the Hilbert spaces $H=\R^d$ and 
$V = \R^{2 \times d}$. For $A \in \R^{d \times D}$ 
 recall also the encoders $\alpha_A:\hat{H}\rightarrow \R^D$ and 
 $\beta_A :\hat{V} \to \R^{2 \times D}$ given respectively by $\alpha_A(x)=(|\ip{x}{a_k}|)_{k\in [D]}$, and  $\beta_A(X)\!=\downarrow\!(XA)
$. 
Our main result reads as follows. 
\begin{theor}\label{th1}
In the case $n=2$, the following are equivalent.
\begin{enumerate}
    \item $\alpha_A$ is injective, hence the columns of $A$ form a phase retrievable frame;
    \item $\beta_A$ is injective, hence $A$ is a universal key.
\end{enumerate}
\end{theor}

\begin{rem}
Perhaps it is not surprising that, if an equivalence between the phase retrieval problem and permutation invariant representations is possible, then this should occur for $n=2$. This statement is suggested by the observation that $O(1)$ is isomorphic with $S_2$, the group of the $2\times 2$ permutation matrices. What is surprising that, in fact, the two embeddings are intimately related, as the proof and corollaries show. 
\end{rem}

\begin{proof}[Proof of Theorem \ref{th1}]

Let $X\in V=\R^{2\times d}$. Denote by $x_1,x_2\in\R^d$ its two rows transposed, that is
\[ 
X = \left[ 
\begin{array}{c}
     x_1^T  \\
     x_2^T
\end{array}
\right].
\]
Notice that, for each $k\in [D]$, the $k^{th}$ column of $\beta_A(X)$ is given by
\begin{align*}
\downarrow (Xa_k)=\begin{bmatrix}
max(\ip{x_1}{a_k},\ip{x_2}{a_k})  \\
min(\ip{x_1}{a_k}, \ip{x_2}{a_k})
\end{bmatrix}.
\end{align*}
The {\em key observations} are the following relationships between $min$, $max$, and the absolute value $|\cdot|$: 
\begin{eqnarray}
\nonumber
   |u-v| &  = & max(u,v) - min(u,v)\\
\nonumber   u+v &  = & max(u,v) + min(u,v) \\
\nonumber   max(u,v) &  = & \frac{1}{2} (u+v+|u-v|)\\
\nonumber    min(u,v) &  = & \frac{1}{2} (u+v-|u-v|) \\
\nonumber |\,|u|-|v|\,| & = & min(|u-v|,|u+v|)
\end{eqnarray}
In particular, these show that:
\begin{align*}
&\begin{bmatrix}
1 & -1 \\
1 &  1
\end{bmatrix}\beta_A(X)=	\begin{bmatrix}
1 & -1 \\
1 &  1
\end{bmatrix} \cdot \downarrow (XA)=\\
&=\begin{bmatrix}
|\ip{x_1-x_2}{a_1}|, & \dots & ,|\ip{x_1-x2}{a_D}|\\
\ip{x_1+x_2}{a_1} ,& \dots &, \ip{x_1-x_2}{a_D}
\end{bmatrix}\\
&= \begin{bmatrix}
(\alpha_A(x_1-x_2))^T\\
(T_A(x_1+x_2))^T
\end{bmatrix}
\end{align*}
Where, $T_A$ was introduced in equation (\ref{eq1}). 
\medskip

$(1) \to (2):$ Suppose that $\alpha_A$ is injective. Let $X=\begin{bmatrix} x_1^T
\\ x_2^T\end{bmatrix}$ and  
       $Y=\begin{bmatrix} y_1^T\\ y_2^T \end{bmatrix}$, such that $\beta_A(X)=\beta_A(Y)$.
Then 
\begin{align*}
&\begin{bmatrix}
1 & -1 \\
1 &  1
\end{bmatrix}\beta_A(X)=	\begin{bmatrix}
1 & -1 \\
1 &  1
\end{bmatrix} \beta_A(Y) \\ \implies& \begin{bmatrix}
    (\alpha_A(x_1-x_2))^T \\
    T_A(x_1+x_2)^T
\end{bmatrix} = \begin{bmatrix}
    (\alpha_A(y_1-y_2))^T \\
    T_A(y_1+y_2)^T
\end{bmatrix}.
\end{align*}
But now,
$\alpha_A(x_1-x_2)=\alpha_A(y_1-y_2)) \implies x_1-x_2=y_1-y_2\text{ or }  x_1-x_2=y_2-y_1$
and
\[ T_A(x_1+x_2)=T_A(y_1+y_2)^T \implies x_1+x_2=y_1+y_2\]
Thus we have that 
\[
\left\{
\begin{array}{rcl}
	x_1 & = & y_1\\
	x_2 & = & y_2
\end{array}
\right\}
\text{ or }
\left\{
\begin{array}{rcl}
	x_1 & = & y_2\\
	x_2 & = & y_1
\end{array}
\right\}
\]

Either case means
\begin{align*}
&\iff X=Y \text{ or } X= \begin{bmatrix}
0 & 1 \\
1 &  0
\end{bmatrix}Y \\ &\iff [X] =[Y]
\end{align*}
So, $\beta_A$ is injective.

$(2) \to (1):$ Suppose that $\beta_A$ is injective. Let $x,y \in \R^d$ such that $\alpha_A(x)=\alpha_A(y)$, i.e. $|\ip{x}{a_k}|=|\ip{y}{a_k}|$, $\forall k \in [D]$.
Let $X=\begin{bmatrix} x^T
\\ -x^T\end{bmatrix}$ and  
       $Y=\begin{bmatrix} y^T\\ -y^T \end{bmatrix}$. Then,
 \[\begin{bmatrix}
1 & -1 \\
1 &  1
\end{bmatrix}\beta_A(X)=\begin{bmatrix}
    \alpha_A(2x)^T \\
    T_A(0)^T
\end{bmatrix}   = 2 \begin{bmatrix}
    \alpha_A(2x)^T \\
    0
\end{bmatrix}\]
and
\[\begin{bmatrix}
1 & -1 \\
1 &  1
\end{bmatrix}\beta_A(X)=\begin{bmatrix}
    \alpha_A(2y)^T \\
    T_A(0)^T
\end{bmatrix}   = 2 \begin{bmatrix}
    \alpha_A(2y)^T \\
    0
\end{bmatrix}\]
Thus  $\beta_A(X)= \beta_A(Y)$. Since $\beta_A$ is assumed injective, it follows that $X=Y$ or $X= \begin{bmatrix}
0 & 1 \\
1 &  0
\end{bmatrix}Y$. So, $x=y$ or $x=-y$. We conclude that $[x]=[y]$, so $\alpha_A$ is injective. 
\medskip

\end{proof}
\begin{cor}
If $\beta_A$ is injective, then $D \geq 2d-1$. 
\end{cor}
\begin{cor}
    If $D=2d-1$, then $\beta_A$ is injective if and only if $A$ is a full spark frame.
\end{cor}
Both results follow necessary and sufficient conditions established in, e.g. \cite{balan06}.
Recall that a frame in $\R^d$ is said {\em full spark} if any subset of $d$ vectors is linearly independent (hence basis). 
\begin{rem}
 Assume $D=2d-1$.  Note the embedding dimension for $\hat{V}=\widehat{\R^{2\times d}}$ is $m=2(2d-1)=4d-2=2\,dim(V)-2$. In particular this shows the minimal dimension of bi-Lipschitz Euclidean embeddings may be smaller than twice the intrinsic dimension of the Hilbert space where the group acts on. Both papers \cite{balan2022permutation} and \cite{cahill22} present (bi)Lipschitz embeddings into 
 $\R^{2\, dim(V)}$.
\end{rem}

\begin{rem}
    As was derived in the proof, $\alpha_A$,  $\beta_A$ and $T_A$ are intimately related:
\begin{equation}\label{eqrem}
\beta_A
\left(
\left[
\begin{array}{c}
x_1^T \\
x_2^T
\end{array}
\right] 
\right)
= \frac{1}{2}
\left[ \begin{array}{cc}
  1   &  1 \\
   -1  &  1
\end{array}
\right]~
\left[
\begin{array}{c}
\alpha_A(x_1-x_2)^T \\
T_A(x_1+x_2)^T
\end{array}
\right]
\end{equation}
In particular, any algorithm for solving the phase retrieval problem solves also the inversion problem for $\beta_A$. Let $\omega_A:\R^{D}\rightarrow \R^d$
denote a left inverse of $\alpha_A$ on the metric space $\widehat{\R^d}$. This means $\omega_A(\alpha_A(x))\sim x$ in $\R^d/O(1)$. 
Denote by $T_A^{\dag}$ a left inverse of the analysis operator (e.g., the synthesis operator associated to the canonical dual frame). Thus $T_A^{\dag}T_A=I_d$. 
Then an inverse for $\beta_A$ is:
\begin{equation}\label{eq:inv}
\beta_A^{-1}(Y) = \frac{1}{2}
\left[ 
\begin{array}{c}
T_A^\dag (y_2) + \omega_A(y_1)  \\
T_A^\dag (y_2) - \omega_A(y_1)  
\end{array}
\right]
\end{equation}
where $Y=\left[ \begin{array}{c} y_1^T \\ y_2^T \end{array} \right]$.
\end{rem}

\begin{rem}
Equations (\ref{eq:inv}) suggest a lower dimensional embedding than $\beta_A$. Specifically, first we compute the average $y_1 = \frac{1}{2}(x_1+x_2)$ which is of size $\R^d$, and then encode the difference $x_1-x_2$ using $\alpha_A$, $y_2=\alpha_A(x_1-x_2)$. We obtain the following modified encoder, $\tilde{\beta}_A:\R^{2\times d}\rightarrow\R^{d+D}$:
\begin{equation}
\tilde{\beta}_A(x)=\left[ 
    \begin{array}{cc} 
    \frac{1}{2}(x_1+x_2)^T & \alpha_A(x_1-x_2)^T
    \end{array}
    \right].
\end{equation}
With the $\omega_A$ left inverse of $\alpha_A$, the inverse of $\tilde{\beta}_A$ is given by:
\begin{equation}\label{eq:inv2} \tilde{\beta}_A^{-1}(Y) = 
\left[ 
\begin{array}{c}
y_1 + \frac{1}{2}\omega_A(y_2)  \\
y_1 - \frac{1}{2}\omega_A(y_2)  
\end{array}
\right]
\end{equation}
where $y_1=Y(1:d)$ and $y_2=Y(d+1:d+D)$. 
In the case when $D=D_{min}=2d-1$, the minimal embedding dimension is $m=d+D=3d-1$ (instead of $4d-2$ or $4d=2\,dim(V)$).
\end{rem}
Reference \cite{balan2022permutation} shows that an upper Lipschitz bound for embedding $\beta_A$ is $\sigma_1(A)$, where
$\sigma_1(A)$ is the largest singular value of $A$. Same reference shows that if 
 $\beta_A$ is injective then there is also a strictly positive lower Lipschitz bound without providing a formula. 
Using \Cref{eqrem} we provide explicit estimates of these bounds.

\begin{theor}\label{lowerL}
Assume $A\in\R^{d\times D}$ is a universal key for $\R^{2\times d}$ (i.e., $\beta_A:\widehat{\R^{2\times d}}\rightarrow\R^{2\times D}$ is injective), or, equivalently (according to Theorem \ref{th1}), 
the columns of $A$ form a phase retrievable frame in $\R^d$ (i.e., $\alpha_A : \widehat{\R^d} \to \R^D$ is injective). Then both $\alpha_A$ and $\beta_A$ are bi-Lipschitz with same Lipschitz constants, where distances are given by $\dispr([x],[y])= \min(\norm{x-y},\norm{x+y})$ on $\hat{H}$, and $\disper([X],[Y])= \min_{P \in S_2} \norm{X-PY}$ on $\hat{V}$, respectively.
The optimal lower and upper Lipschitz constants are given by: 
\begin{equation}
    \label{eq:A0}
A_0 = \min_{I\subset[D]} \sqrt{\sigma_d^2(A[I]) + \sigma_d^2(A[I^c])} ~~,~~B_0=\sigma_1(A)
\end{equation}
where $\sigma_1(A)$ is the largest singular value of $A$ (equals the square-root of upper frame bound) and $\sigma_d(A[J])$ is the $d^{th}$ singular value of submatrix of $A$ indexed by $J$. Furthermore, these bounds are achieved by the following vectors. Let $I_0$ denote a optimal partition in (\ref{eq:A0}) and let $u_1$, $u_2$  denote the normalized left singular vectors of $A[I_0]$ and $A[I_0^c]$, respectively, each associated to the $d^{th}$ singular value. Let $u$ be the normalized principal left singular vector associated to $A$.(i.e., associated to the largest singular value). Then:

\begin{enumerate}
    \item The upper Lipschitz constant $B_0$ is achieved as follows: (i) for map $\alpha_A$ by vectors $x_{\max}=u$ and $y_{\max}=0$; (ii) for map $\beta_A$ by vectors 
    $X_{\max}=\left[ \begin{array}{c}
    u^T \\ 0 
    \end{array} \right]$ and
    $Y_{\max}=0$.
    \item The lower Lipschitz constant $A_0$ is achieved as follows: (i) for map $\alpha_A$ by vectors $x_{\min}=u_1+u_2$ and $y_{\min}=u_1-u_2$; (ii) for map $\beta_A$ by vectors 
    $X_{\min}=\left[ \begin{array}{c} 
    (u_1 + u_2)^T \\
    0 
    \end{array} \right]$ and
    $Y_{\min}=\left[ \begin{array}{c}  
    u_1^T \\
    u_2^T
    \end{array} \right]$.
\end{enumerate}
\end{theor}

\begin{rem} 
The optimal Lipschitz constants for the map $\alpha_A$ were obtained in \cite{balwan,balazou}, including the optimizers. However, for reader's convenience, we prefer to give direct proofs of these results.  
\end{rem}

\begin{proof}
\begin{enumerate}
    \item Upper Lipschitz constants. 
    
    (i) Let $x,y \in \R^{d}$. Then
\begin{align*}
 & \scriptstyle \norm{\alpha_A(x)-\alpha_A(y)}^2= \sum_{i=1}^D||\ip{a_i}{x}|-|\ip{a_i}{y}||^2 = \\  & = \scriptstyle \sum_{i=1}^D \min\left(|\ip{a_i}{x-y}|^2,|\ip{a_i}{x+y}|^2\right)
 \le \\
 & \scriptstyle \le \min\left(\sum_{i=1}^D |\ip{a_i}{x-y}|^2,\sum_{i=1}^D |\ip{a_i}{x+y}|^2 \right)
 \\ & \scriptstyle \le \sigma_1^2(A) \dispr([x],[y])^2.
\end{align*}
So $\sigma_1(A)$ is an upper Lipschitz bound for the map $\alpha_A$.
Now for $x_{\max}=u,y_{\max}=0$ notice that
\begin{align*}
 &  \scriptstyle \norm{\alpha(x_{\max})-\alpha(y_{\max})}^2=  \sum_{i=1}^D|\ip{a_i}{u}|^2 =\\ & \scriptstyle =\sigma_1^2(A) \norm{u}^2 =\sigma_1^2(A) \dispr([x_{\max}],[y_{\max}])^2.
\end{align*}
Thus, the upper Lipschitz constant $\sigma_1(A)$ is in fact optimal (tight).

(ii) Map $\beta_A$. Let $X,Y \in \R^{2 \times D}$ and $P_0 \in S_2$ be a permutation that achieves the distance between $X$ and $Y$, i.e. $\norm{X-P_0Y}= \disper([X],[Y])$. Note that
 \begin{align*}
& \scriptstyle\norm{\beta_{A}(X)-\beta_{A}(Y)}^2 = \sum_{k=1}^D \norm{(\Pi_k X -\Xi_kY)a_k}^2= \\  & \scriptstyle =\sum_{k=1}^D \norm{(\Xi_k^T\Pi_k X -Y)a_k}^2
  \end{align*}
  for some $\Pi_k,\Xi_k \in S_2$ that align the vectors.
From rearrangement lemma we have that \[\scriptstyle \norm{(\Pi_k X -\Xi_kY)a_k}\leq \norm{(X-P_0Y)a_k},~\forall k \in [D]\] 
so, 
\begin{align*}
 \scriptstyle \sum_{k=1}^D \norm{(\Xi_k^T\Pi_k X -Y)a_k}^2 \le& \norm{A}^2\norm{\norm{X-P_0Y}}^2\\ =& \scriptstyle \sigma_1^2(A)  \disper([X],[Y])^2.
\end{align*}
Therefore, we conclude that $\sigma_1(A)$ is an upper Lipschitz constant for map $\beta_A$.
We still need to show that this bound is achieved (i.e., it is optimal).
For $X_{\max}$ and $Y_{\max}$ defined in part 1) of \cref{lowerL},
    \begin{align*}
& \scriptstyle \norm{\beta_{A}(X_{\max})-\beta_{A}(Y_{\max})}^2= \norm{\beta_{A}(X_{\max})}^2 =\sum_{k=1}^D \ip{u}{a_k}^2  =  \sigma_1^2(A).
  \end{align*}
  and $\disper(X_{\max},Y_{\max})=1$. 
Thus $B_0$ is the optimal Lipschitz constant both for $\alpha_A$ and for $\beta_A$.

\item Lower Lipschitz constants. 

(i) Let $x,y\in\R^d$ and define the auxiliary set
\[\scriptstyle S = S(x,y) := \{j\in [D]~:|~\ip{x-y}{a_j}| \le |\ip{x+y}{a_j}| \} \]
Then
\begin{align*}
   & \scriptstyle \norm{\alpha(x)-\alpha(y)}^2= \sum_{i=1}^D\left| |\ip{a_i}{x}|-|\ip{a_i}{y}| \right|^2 = \\ & = \scriptstyle \sum_{i \in S}|\ip{a_i}{x-y}|^2+\sum_{i \in S^c}|\ip{a_i}{x+y}|^2 \geq \\ & \scriptstyle\sigma_d^2(A[S])+\sigma_d^2(A[S^c]) \dispr([x],[y])^2 \geq A_0^2 \dispr([x],[y])^2 .
\end{align*}
So $A_0$ is a lower Lipschitz bound for $\alpha_A$, but we still need to show that it is optimal.

Let $I_0$ be the optimal partition, and let $u_1$, $u_2$ be normalized left singular vectors as in the statement of \Cref{lowerL}. Then:

\begin{align*}
 & \scriptstyle \norm{\alpha_A(u_1+u_2)-\alpha_A(u_1-u_2)}^2= \\ & \scriptstyle = \sum_{i=1}^D||\ip{a_i}{u_1+u_2}|-|\ip{a_i}{u_1-u_2}||^2 = \\  & = \scriptstyle \sum_{i=1}^D \min\left( |\ip{a_i}{2u_2}|^2,|\ip{a_i}{2u_1}|^2 \right)
 \le \\
 & \scriptstyle \le 4\,\left( \sum_{i\in I_0} |\ip{a_i}{u_1}|^2 + \sum_{i\in I_0^c} |\ip{a_i}{u_2}|^2 \right)
 \\ & \scriptstyle =4 (\sigma_d^2(A[I_0]) + \sigma_d^2(A[I_0^c])) = A_0^2 \dispr([u_1+u_2],[u_1-u_2])^2,
\end{align*}
where we used again that $|\,|a|-|b|\,|=min(|a-b|,|a+b|)$ for any two real numbers $a,b\in\R$, and, for the inequality, at every $i\in[D]$ we made a choice between the two terms. Since the reverse inequality is also true, it follows that $x_{\min}=u_1+u_2$ and $y_{\min}=u_1-u_2$ achieve the lower bound $A_0$ for $\alpha_A$.

(ii) Consider now the map $\beta_A$. Let $X,Y\in\R^{2\times d}$ and define the auxiliary set
\begin{align*}
\scriptstyle S = S(X,Y) := \{& \scriptstyle j\in [D]~:|~\ip{x_1-x_2-y_1+y_2}{a_j}| \le 
\\ & \scriptstyle \le |\ip{x_1-x_2+y_1-y_2}{a_j}| \}
\end{align*}

Then, using \Cref{eqrem} we have that 
\[\scriptstyle
\norm{\beta_A(X)-\beta_A(Y)}^2= \]
\[ \scriptstyle \hspace{-3mm} \frac{1}{2}\left(\norm{\alpha_A(x_1-x_2)-\alpha_A(y_1-y_2)}^2+\norm{T_A(x_1+x_2-y_1-y_2)}^2\right) = \]
\[ \scriptstyle \hspace{-3mm}
\frac{1}{2}\sum_{j \in S}|\ip{x_1-x_2 - y_1+y_2}{a_j}|^2 +|\ip{x_1+x_2 - y_1-y_2}{a_j}|^2 + \]
\[ \scriptstyle \hspace{-5mm}
+\frac{1}{2}\sum_{j\in S^c} |\ip{x_1-x_2 + y_1-y_2}{a_j}|^2+ |\ip{x_1+x_2 -y_1-y_2}{a_j}|^2 =\] \[  \scriptstyle
= \sum_{j\in S} |\ip{x_1-y_1}{a_j}|^2 + |\ip{x_2-y_2}{a_j}|^2 + \]
\[\scriptstyle + \sum_{j\in S^c} |\ip{x_1-y_2}{a_j}|^2 + |\ip{x_2-y_1}{a_j}|^2 \geq \]
\[ \scriptstyle \geq \sigma_d^2(A[S])\left(\norm{x_1-y_1}^2+\norm{x_2-y_2}^2\right) +\]\[ \scriptstyle + \sigma_d^2(A[S^c]) \left(\norm{x_1-y_2}^2+\norm{x_2-y_1}^2\right) \geq \]
\[ \scriptstyle \geq A_0^2 \disper([X],[Y])^2. \]
Therefore $A_0$ is a lower Lipschitz constant for $\beta_A$. 

It remained to prove that this bound is tight, i.e., it is achieved. Let $X_{\min}$ and $Y_{\min}$ be as in the statement of \Cref{lowerL}.
 Then
 \[ \scriptstyle
\norm{\beta_A(X_{\min})-\beta_A(Y_{\min})}^2= \]
\[ \scriptstyle \frac{1}{2}\left(\norm{\alpha_A(u_1+u_2)-\alpha_A(u_1-u_2)}^2+\norm{T_A(u_1+u_2-u_1-u_2)}^2\right) = \]
\[ \scriptstyle \frac{1}{2}\left(\norm{\alpha_A(u_1+u_2)-\alpha_A(u_1-u_2)}^2 \right) = A_0^2 \disper([X_{\min}],[Y_{\min}])^2 \]
 where the last equality follows from the fact that the lower Lipschitz constant of $\alpha_A$ is
achieved by $u_1+u_2$ and $u_1-u_2$, and the fact that
 $\disper([X_{\min}],[Y_{\min}])^2=2$.

 So $A_0$ is indeed the optimal lower Lipschitz constant for $\beta_A$.

\end{enumerate}
\end{proof}

\section{Conclusion}

In this paper we analyzed two representation problems, one arising in the phase retrieval problem
and the other one in the context of permutation invariant representations. We showed that the real phase retrieval problem in a finite dimensional vector space $H$ is entirely equivalent to the permutation invariant representations for the space $V=\R^{2\times dim(H)}$. 
Our analysis proved that phase retrievability is equivalent to the universal key property in the case of encoding ${2\times d}$ matrices.
This result is derived based on the lattice space structure $(\R,+,min,max)$.
It is still an open problem to understand the relationship between $\alpha_A$ and $\beta_A$ in the case $n>2$. 
A related problem is the implementation of the sorting operator using a neural network that has ReLU as activation function (or, even the absolute value $|\cdot|$). 
Efficient implementations of such operator may yield novel relationships between $\alpha_A$ and $\beta_A$, in the case $n\geq 3$.


\section*{Acknowledgment}

The authors have been supported in part by a NSF award under grant DMS-2108900 and
by the Simons Foundation.

\bibliographystyle{IEEEtran.bst}
\bibliography{references.bib}

\end{document}